\documentclass[12pt,oneside]{article}

\usepackage{amssymb}
\usepackage{amsxtra}
\usepackage{amsmath}
\usepackage{amsfonts}
\usepackage{CJK}
\usepackage[titletoc]{appendix}
\usepackage{graphicx}
\usepackage{color}
\numberwithin{equation}{section}

\newcommand{\eqa}{\begin{eqnarray}}
\newcommand{\eeqa}{\end{eqnarray}}
\newcommand{\beq}{\begin{equation}}
\newcommand{\eeq}{\end{equation}}
\newcommand{\nn}{\nonumber}

\allowdisplaybreaks

\begin{document}

\date{}
\author{Fang Fang$^{1}$ , Beibei Hu$^{1,}$\thanks{Corresponding author. E-mails: hu\_chzu@163.com } , Ling Zhang$^{1}$, Ning Zhang$^{2}$\\
\small \textit{1.School of Mathematics and Finance, Chuzhou University, Anhui, 239000, China}\\
\small \textit{2.Department of Basical Courses, Shandong University of Science and Technology, Taian 271019, China}}
\title{Riemann-Hilbert approach for a mixed coupled nonlinear Schr\"{o}dinger system and its soliton solutions
}
\maketitle
 \hrulefill

\begin{abstract}

In this work, we examine the integrable mixed coupled nonlinear Schr\"{o}dinger (mCNLS) system, which describe the propagation of an optical
pulse in a birefringent optical fiber. By the Riemann-Hilbert(RH) approach, the N-soliton solutions of the mCNLS system can be expressed explicitly when the jump matrix of a specific RH problem is a $3\times3$ unit matrix. As a special example, the expression of one- and two-soliton are displayed explicitly. More generally, as a promotion, an integrable generalized multi-component NLS system with its linear spectral problem be discussed. It is hoped that our results can help enrich the nonlinear dynamical behaviors of the mCNLS.
\\
\\
\textbf{Keywords:}~~Riemann-Hilbert approach; mixed coupled nonlinear Schr\"{o}dinger (mCNLS) system; soliton solution; boundary conditions.
\end{abstract}
\hrulefill

\section{Introduction}

It is well known that solving nonlinear evolution equations becomes a challenging task due to the complexity of nonlinear systems. In particular, the acquisition of precise solutions is crucial for research in various fields. Through the many years of efforts of mathematics and physicists, a variety of construction methods for precise solutions have been established, such as the inverse scattering method \cite{GGKM1967,AKNS1981,RB1984,PAC1989,AC2010,GB2015,JJL2017}, the Hirota's bilinear method \cite{HR1971,ZY2017,HZ2017,MWX20181}, the B\"{a}cklund transformation method \cite{PJO1986}, the Darboux transformation (DT) method \cite{VBM1991} and others \cite{LSY1990,LJB2000,FEG2003}.
In recent years, with the development of the soliton theory, more and more scholars are paying attention to Riemann-Hilbert (RH) method\cite{YJK2010}, which is a new powerful method for solving integrable linear and nonlinear partial differential equations (PDEs) \cite{GBL2012,WDS2014, ZYS2017,GXG2016,WZ2016,WJP20171,MWX2018,HJ2018,ALP2018}.  The main idea of the this method is to established a corresponding matrix RH problem on the Lax pair of integrable equations. Furthermore, the RH method is also an effective way to examined the initial boundary value problems \cite{YZY2017,ZN2017,XBQ2018,Hu1,Hu2,Hu3,Hu4} and the long-time asymptotic behavior \cite{DP1993} of the integrable nonlinear evolution PDEs.

The NLS equation is one of the most paramount integrable systems in mathematics and physics which reads
\beq iq_{t}\pm q_{xx}+|q|^2q=0, \label{1.1}\eeq
and which arises in various physical backgrounds involving fluid mechanics, nonlinear optics, plasma physics, Bose-Einstein condensation and other fields
whose describes wave dynamics of nonlinear pulses propagation under the absence of optical losses in monomode fiber.
However, a slice of phenomena have been observed by experimental which cannot be explained by Eq.\eqref{1.1}.
In order to explain this phenomena, quiet a few researchers examined the the two-component case (known as the Manakov system)
\beq\left\{\begin{array}{l}
iq_{1t}+\frac{1}{2}q_{1xx}+\epsilon(|q_1|^2+|q_2|^2)q_1=0,\\
iq_{2t}+\frac{1}{2}q_{2xx}+\epsilon(|q_1|^2+|q_2|^2)q_2=0,
\end{array}\right.\epsilon=\pm1.\label{1.2}\eeq
Here, $\epsilon=1$ and $\epsilon=-1$ means the defocusing case and the focusing case, respectively. This system was
first introduced by Manakov \cite{MSV1974} to describe the propagation of an optical pulse in a birefringent
optical fibre. The Manakov system \eqref{1.2} affords the mathematical luxury of extending the local
linearized analysis to a construction of the entire nonlinear unstable manifold of the underlying
oscillatory wave.

In the present paper, for convenience and based on RH method, we consider the following coupled focusing-defocusing NLS system:
\beq\left\{\begin{array}{l}
iq_{1t}+\frac{1}{2}q_{1xx}+(|q_2|^2-|q_1|^2)q_1=0,\\
iq_{2t}+\frac{1}{2}q_{2xx}+(|q_2|^2-|q_1|^2)q_2=0.
\end{array}\right.\label{1.3}\eeq
where $q_1$ and $q_2$ are two-component electric field functions; $x$ and $t$ denote the direction of propagation and time variables, respectively. Generally, system \eqref{1.3} is called the mixed coupled nonlinear Schr\"{o}dinger (mCNLS) equations, it is a completely integrable equations, and quiet a few properties of the system \eqref{1.3} have been discussed. As an example, soliton collisions with a
shape change have been investigated by intensity redistribution in \cite{KT2006}. The bright-dark solitons
and their collisions have been studied in mixed N-coupled NLS equations in
\cite{VM2008}. Very recently, vector rogue wave (type I and type II) solutions and bright-dark rogue wave
solutions \cite{LLM2016} have been given by using the DT. The initial-boundary value problems have been studied by the Fokas method \cite{TSF2016}. However, to the best of the author's knowledge, the N-soliton solution of system \eqref{1.3} via the Riemann-Hilbert approach have not been investigated before.

The organization of this paper is as follows. In section 2, we will construct a specific RH problem based on the inverse scattering transformation. In section 3, we compute N-soliton solutions of the mCNLS system from a specific RH problem,
which possesses the identity jump matrix on the real axis. In section 4, as a promotion, We briefly explained an integrable generalized multi-component NLS system with its linear spectral problem can be solved by the same ways. And some conclusions
are given in the final section.

\section{ The Riemann-Hilbert problem}

Consider the following Lax pair of system \eqref{1.3}
\begin{subequations}
\begin{align}
&\Phi_x=M(x,t,k)\Phi=(-ik\Lambda+iQ)\Phi,\label{2.1a}\\
&\Phi_t=N(x,t,k)\Phi=(-ik^2\Lambda+ikQ+\frac{1}{2}(iQ^2-Q_x))\Phi,\label{2.1b}
\end{align}
\end{subequations}
where
\eqa\begin{array}{l}
\Lambda=\left(\begin{array}{ccc}
-1&0&0\\
0&1&0\\
0&0&1\end{array} \right),
Q=\left(\begin{array}{ccc}
0&-q_1^*&q_2^*\\
q_1&0&0\\
q_2&0&0\end{array} \right).
\end{array}\label{2.2}\eeqa

It is easy to see that \eqref{2.1a} and \eqref{2.1b} is equivalent to
\begin{subequations}
\begin{align}
&\Phi_x+ik\Lambda\Phi=Q_1\Phi,\label{2.3a}\\
&\Phi_t+ik^2\Lambda\Phi=Q_2\Phi,\label{2.3b}
\end{align}
\end{subequations}
where
\beq Q_1=iQ,\,\,Q_2=ikQ+\frac{1}{2}(iQ^2-Q_x).\label{2.4}\eeq

Obviously, $\tilde{A}(x,t,k)=e^{-(ik\Lambda x+ik^2)\Lambda t}$ is a solution for the Lax pair \eqref{2.3a} and \eqref{2.3b} at this time. Let $\Psi(x,t,k)=J(x,t,k)\tilde{A}(x,t,k)$ , then the spectral problems about $J(x,t,k)$ are defined as
\begin{subequations}
\begin{align}
&J_x+ik[\Lambda,J]=Q_1J,\label{2.5a}\\
&J_t+ik^2[\Lambda,J]=Q_2J.\label{2.5b}
\end{align}
\end{subequations}

Now, we are construct two Jost solutions $J_{\pm}= J_{\pm}(x,k)$ of Eq.\eqref{2.5a} for $k\in \mathbb{R}$
\eqa J_{+}=([J_{+}]_1,[J_{+}]_2,[J_{+}]_3),\label{2.6}\eeqa
\eqa J_{-}=([J_{-}]_1,[J_{-}]_2,[J_{-}]_3),\label{2.7}\eeqa
with the boundary conditions
\begin{subequations}
\begin{align}
&J_{+}\rightarrow \mathrm{I},\,\,x\rightarrow -\infty,\label{2.8a}\\
&J_{-}\rightarrow \mathrm{I},\,\,x\rightarrow +\infty.\label{2.8b}
\end{align}
\end{subequations}
Where $[J_{\pm}]_n(n=1,2,3)$ denote the $n$-th column vector of $J_{\pm}$, $\mathrm{I}=diag\{1,1,1\}$ is a $3\times3$ identity matrix,
and the subscripts of $J(x,k)$ represent which end of the $x$-axis the boundary
conditions are set. Infact, the two Jost solutions $J_{\pm}= J_{\pm}(x,k)$ of Eq.\eqref{2.5a} for $k\in \mathbb{R}$ are determined by the following Volterra integral equations
\eqa J_{+}(x,k)=\mathrm{I}-\int_x^{+\infty}e^{-ik\hat\Lambda (x-\xi)}Q_1(\xi)J_{+}(\xi,k)d\xi,\label{2.9}\eeqa
\eqa J_{-}(x,k)=\mathrm{I}+\int_{-\infty}^xe^{-ik\hat\Lambda (x-\xi)}Q_1(\xi)J_{-}(\xi,k)d\xi,\label{2.10}\eeqa
where $\hat\Lambda$ represents a matrix operator acting on $3\times3$ matrix $X$ by $\hat\Lambda X=[\Lambda,X]$ and by
$e^{x\hat\Lambda}X=e^{x\Lambda}Xe^{-x\Lambda}$.

Moreover, after simple analysis, we find that $[J_{+}]_1, [J_{-}]_2$ and $[J_{-}]_3$ admits analytic extensions to $C_{+}$. On the other hand, $[J_{-}]_1, [J_{+}]_2$ and $[J_{+}]_3$ admits analytic extensions to the $C_{-}$,
here $C_{+}$ and $C_{-}$ denote the upper half $k$-plane and the lower half $k$-plane, respectively.

Next, let us investigate the properties of ${J}_{\pm}$. Due to the Abel's identity and $\mathrm{Tr}(Q)=0$, the determinants of ${J}_{\pm}$ are constants for all $x$. From the boundary conditions Eq.(2.8), we have
\beq \det {{J}_{\pm }}=1,\quad k \in \mathbb{R}. \label{2.11}\eeq
Introducing a new function $ A(x,k)=e^{-ik\Lambda x}$, we find that spectral problem Eq.\eqref{2.5a} exists two fundamental matrix solutions ${{J}_{+}}A$ and ${{J}_{-}}A$, which are not independent and are linearly associated by a $3\times3$ scattering matrix $S(k)$
\beq {{J}_{-}}A={{J}_{+}}A\cdot S(k),\quad k \in \mathbb{R}.\label{2.12}\eeq
or
\beq {{J}_{-}}={{J}_{+}}A\cdot S(k)A^{-1},\quad k \in \mathbb{R}.\label{2.13}\eeq
It follows from Eq.\eqref{2.11} and \eqref{2.12} we know that
\beq \det S(k)=1. \label{2.14}\eeq
Moreover, let $x$ go to $+\infty$, the $3\times3$ scattering matrix $S(k)$ is given as
\beq  S(k)={{({{s}_{ij}})}_{3\times 3}}=\lim_{x\rightarrow+\infty}A^{-1}J_{-}A=\mathrm{I}+\int_{-\infty}^{+\infty}e^{ik\hat\Lambda \xi}Q_1J_{-}d\xi,\,\,k \in \mathbb{R}. \label{2.15}\eeq
From the analytic property of $J_{-}$, we find that $s_{22},s_{23},s_{32}$ and $s_{33}$ can be analytically extended to $C_{+}$, $s_{11}$
allow analytic extensions to $C_{-}$. Generally speaking, $s_{12},s_{13},s_{21}$ and $s_{31}$ cannot be extended off the
real $x$-axis.

In order to obtain behavior of Jost solutions for very large $k$, we substituting the following expansion
\beq J=J_0+\frac{J_1}{k}+\frac{J_2}{k^2}+\frac{J_3}{k^3}+\frac{J_4}{k^4}+\cdots \quad k\rightarrow\infty,\label{2.16}\eeq
into the Eq.\eqref{2.5a} and comparing the coefficients of the same order of $k$ yields
\begin{subequations}
\begin{align}
& O(k^1):i[\Lambda,J_0]=0, \label{2.17a} \\
& O(k^0):J_{0,x}+i[\Lambda,J_1]-Q_1J_0=0,  \label{2.17b}\\
& O(k^{-1}):J_{1,x}+i[\Lambda,J_2]-Q_1J_1=0,\label{2.17c}
\end{align}
\end{subequations}
From $O(k^1)$ and $O(k^0)$ we have
\beq i[\Lambda,J_1]=Q_1J_0,\,\,J_{0,x}=0. \label{2.18}\eeq

In order to construct the RH problem of the mCNLS system, we must to define another new Jost solution for Eq.\eqref{2.5a} by
\beq P_{+}=([J_{+}]_1,[J_{-}]_2,[J_{-}]_3)=J_{+}AS_{+}A^{-1}=J_{+}A\left(\begin{array}{ccc}
1 & s_{12} & s_{13}\\
0 & s_{22} & s_{23}\\
0 & s_{32} & s_{33}
\end{array}\right)A^{-1}, \label{2.19}\eeq
which is analytic for $k\in C_{+}$ and admits asymptotic behavior for very large $k$ as
\beq P_{+}\rightarrow\mathrm{I},\, k\rightarrow +\infty,\, k\in C_{+}. \label{2.20}\eeq

Furthermore, to obtain the analytic counterpart of $P_{+}$ in $C_{-}$, denoted by$P_{-}$, we consider the adjoint scattering
equation of Eq.\eqref{2.5a}:
\beq Z_x+ik[\Lambda,Z]=-ZQ_1. \label{2.21}\eeq
Obviously, the inverse matrices $J_{\pm}^{-1}$ defined as
\beq {[J_{+}]}^{-1}=\left(\begin{array}{ccc}
{[J_{+}^{-1}]}^1\\
{[J_{+}^{-1}]}^2\\
{[J_{+}^{-1}]}^3
\end{array}\right),\quad
{[J_{-}]}^{-1}=\left(\begin{array}{ccc}
{[J_{-}^{-1}]}^1\\
{[J_{-}^{-1}]}^2\\
{[J_{-}^{-1}]}^3
\end{array}\right),\label{2.22}\eeq
satisfy this adjoint equation \eqref{2.21}, here $[J_{\pm}^{-1}]^n(n=1,2,3)$ denote the $n$-th row vector of $J_{\pm}^{-1}$. Then we can see that ${[J_{+}^{-1}]}^1, {[J_{-}^{-1}]}^2$ and ${[J_{-}^{-1}]}^3$ admits analytic extensions to $C_{-}$. On the other hand, ${[J_{-}^{-1}]}^1, {[J_{+}^{-1}]}^2$ and ${[J_{+}^{-1}]}^3$ admits analytic extensions to the $C_{+}$.

In addition, it is not difficult to find that the inverse matrices $J_{+}^{-1}$ and $J_{-}^{-1}$ satisfy the following boundary conditions .
\begin{subequations}
\begin{align}
&J_{+}^{-1}\rightarrow \mathrm{I},\,\,x\rightarrow -\infty,\label{2.23a}\\
&J_{-}^{-1}\rightarrow \mathrm{I},\,\,x\rightarrow +\infty.\label{2.23b}
\end{align}
\end{subequations}
Therefore, one can define a matrix function $P_{-}$ is expressed as follows:
\beq P_{-}=\left(\begin{array}{ccc}
{[J_{+}^{-1}]}^1\\
{[J_{-}^{-1}]}^2\\
{[J_{-}^{-1}]}^3
\end{array}\right).\label{2.24}\eeq
By techniques similar to those used above, one can show that the adjoint Jost solutions $P_{-}$ are analytic in $C_{-}$ and
\beq P_{-}\rightarrow\mathrm{I},\, k\rightarrow -\infty,\, k\in C_{-}. \label{2.25}\eeq
Assume that $R(k)=S^{-1}(k)$, we have
\beq {J_{-}^{-1}}=AR(k)A^{-1}{J_{+}^{-1}},\label{2.26}\eeq
and
\beq P_{-}=\left(\begin{array}{ccc}
{[J_{+}^{-1}]}^1\\
{[J_{-}^{-1}]}^2\\
{[J_{-}^{-1}]}^3
\end{array}\right)=AR_{+}A^{-1}J_{+}^{-1}=A\left(\begin{array}{ccc}
1 & 0 &0\\
r_{21} & r_{22} & r_{23}\\
r_{31} & r_{32} & r_{33}
\end{array}\right)A^{-1}J_{+}^{-1}, \label{2.27}\eeq

Hence we have constructed two matrix functions $P_{+}(x,k)$ and $P_{-}(x,k)$ which are analytic for $k$ in $C_{+}$ and $C_{-}$, respectively.
In fact, these two matrix functions $P_{+}(x,k)$ and $P_{-}(x,k)$ which can be construct a RH problem:
\beq P_{-}(x,k)P_{+}(x,k)=T(x,k),\,\, k\in C_{-}. \label{2.28}\eeq
where
\beq T(x,k)=AR_{+}S_{+}A^{-1}=\left(\begin{array}{ccc}
1 & s_{12}e^{2ikx} & s_{13}e^{2ikx}\\
r_{21}e^{-2ikx} & 1 & 0\\
r_{31}e^{-2ikx} & 0 & 1
\end{array}\right),\,\, k\in C_{-}. \label{2.29}\eeq
Here we have adopted the identity ${{r}_{11}}{{s}_{11}}+{{r}_{12}}{{s}_{21}}+{{r}_{13}}{{s}_{31}}=1$, and the jump contour is real $x$-axis.

Furthermore, since $J_{-}$ satisfies the temporal part of spectral equation
\eqa J_{-,t}+ik^2[\Lambda,J_{-}]=Q_2J_{-}, \label{2.30}\eeqa
we have
\eqa J_{-}A=J_{+}AS,\,\,(J_{+}AS)_t+ik^2[\Lambda,J_{+}AS]=Q_2J_{+}AS, \label{2.31}\eeqa
suppose $q_1$ and $q_2$ sufficient smoothness and decay as $x\rightarrow\infty$, we have $Q_2\rightarrow 0$ as $x\rightarrow\pm\infty$.
Then taking the limit $x\rightarrow+\infty$ of Eq.\eqref{2.31} yields
\eqa S_{t}=-ik^2[\Lambda,S],\label{2.32}\eeqa
This above equation imply that the scattering data $s_{11},s_{22},s_{33},s_{23},s_{32}$ are time independent, that is to say
\eqa s_{11,t}=s_{22,t}=s_{33,t}=s_{23,t}=s_{32,t}=0,\label{2.33}\eeqa
and the other scattering data satisfies
\eqa && s_{12}(t,k)=s_{12}(0,k)e^{2ik^2t},\,\,s_{13}(t,k)=s_{13}(0,k)e^{2ik^2t}, \nn\\&&
s_{21}(t,k)=s_{21}(0,k)e^{-2ik^2t},\,\,s_{31}(t,k)=s_{31}(0,k)e^{-2ik^2t}.
\label{2.34}\eeqa

\section{The soliton solutions}

In fact, the solution to this RH problem will not be unique unless the zeros of det $P_{+}$ and det $P_{-}$ in the upper and
lower half of the $k$-plane are also specified, and the kernel structures of $P_{\pm}$ at these zeros are provided.
From the definitions of $P_{+}$ and $P_{-}$ as well as the scattering relations between $J_{+}$ and $J_{-}$, we see that
\beq \mathrm{det} P_{+}(x,k)=r_{11}(k),\quad \mathrm{det} P_{-}(x,k)=s_{11}(k),\label{3.1}\eeq
where $r_{11}=s_{22}s_{33}-s_{23}s_{32}$, which imply that the zeros of det $P_{+}$ and det $P_{-}$ are the same as $r_{11}(k)$ and $s_{11}(k)$, respectively. Indeed, owing to the scattering data $s_{11}$ and $r_{11}$ are time independent, then the roots of $s_{11}=0$ and $r_{11}=0$ are also time independent. Furthermore, owing to
$$Q^\dag=\sigma Q\sigma,$$
where $\sigma=diag\{1,-1,1\}$.

It is easy to see that
\beq J_{\pm}^{-1}(x,t,k)=\sigma J_{\pm}^{\dag}(x,t,k^*)\sigma,\label{3.2}\eeq
applying these two reduction conditions to (2.7), then
\beq S^{-1}(k)=\sigma S^{\dag}(k^*)\sigma,\quad P_{-}(x,k)=\sigma P_{+}^{\dag}(x,k^*)\sigma,\label{3.3}\eeq

Suppose that $r_{11}$ possess $N\geq0$ possible zeros in $C_{+}$ denoted by $\{k_m,1\leq m\leq N\}$, and $s_{11}$ possess $N\geq0$ possible zeros in $C_{-}$ denoted by $\{\hat k_m,1\leq m\leq N\}$. For simplicity, we assume that all zeros $\{(k_m,\hat k_m),m=1,2,...,N\}$ are simple zeros of $r_{11}$ and $s_{11}$, which is the generic case. In this case, each of ker $P_{+}(k_m)$ and ker $P_{-}(\hat k_m)$ contains only a single column vector $v_{m}$ and row vector $\hat v_{m}$, respectively, such that
\eqa P_{+}(k_m)v_{m}=0,\quad \hat v_{m}P_{-}(\hat k_m)=0.  \label{3.4}\eeqa

Owing to $P_{+}(k)$ is the solution of spectral problem \eqref{2.5a}, we assume that the asymptotic expansion of $P_{+}(k)$ at large $k$ as
\beq P_{+}=\mathrm{I}+\frac{P_{+}^{(1)}}{k}+O(k^{-2})\quad k\rightarrow\infty,\label{3.5}\eeq
substitute the above expansion into \eqref{2.5a} and \eqref{2.5b} and compare $O(1)$ terms obtain
\beq Q_{1}=i[\Lambda,P_{+}^{(1)}]=\left(\begin{array}{ccc}
0 & -2i(P_{+}^{(1)})_{12} & -2i(P_{+}^{(1)})_{13}\\
2i(P_{+}^{(1)})_{21} & 0 & 0\\
2i(P_{+}^{(1)})_{31} & 0 & 0
\end{array}\right).\label{3.6}\eeq
then the potential functions $q_1$ and $q_2$ can be reconstructed by
\beq q_1=2(P_{+}^{(1)})_{21},\quad q_2=2(P_{+}^{(1)})_{31},\label{3.7}\eeq
where $P_{+}^{(1)}=(P_{+}^{(1)})_{3\times3}$ and $(P_{+}^{(1)})_{ij}$ is the $(i;j)$-entry of $P_{+}^{(1)},i,j=1,2,3$.

In order to obtain the spatial evolutions for vectors $v_m(x,t)$, on the one hand, we taking the $x$-derivative to equation $P_{+}v_m=0$ and using \eqref{2.5a} obtain
\beq P_{+}v_{m,x}+ik_mP_{+}\Lambda v_m=0,\label{3.8}\eeq
thus
\beq v_{m,x}=-ik_m\Lambda v_m,\label{3.9}\eeq
on the other hand, we also taking the $t$-derivative to equation $P_{+}v_m=0$ and using (2.5b) obtain
\beq P_{+}v_{m,t}+ik_m^2P_{+}\Lambda v_m=0,\label{3.10}\eeq
thus
\beq v_{m,t}=-ik_m^2\Lambda v_m,\label{3.11}\eeq
By solving \eqref{3.9} and \eqref{3.11} explicitly, we get
\beq v_{m}(x,t)=e^{-ik_m\Lambda x-ik_m^2\Lambda t}v_{m0},\label{3.12}\eeq
\beq \hat v_{m}(x,t)=v_m^\dag\sigma=\hat v_{m0}e^{ik_m^*\Lambda x+ik_m^{*2}\Lambda t}\sigma,\label{3.13}\eeq
where $\sigma=diag\{1,-1,1\}$.

In order to construct multi-soliton solutions for the mCNLS system \eqref{1.3}, one can choose the jump matrix $T=\mathrm{I}$ is a $3\times3$ unit matrix in \eqref{2.28}. That is to say, the discrete scattering data $r_{12}=r_{13}=s_{21}=s_{31}=0$, consequently, the unique solution to this special RH problem have been solved in \cite{YJK2010}, and the result is
\begin{subequations}
\begin{align}
&P_{+}(k)=\mathrm{I}-\sum_{m=1}^{N}\sum_{n=1}^{N}\frac{v_m(M^{-1})_{mn}\hat v_n}{k-\hat k_m},\label{4.1a}\\
&P_{-}(k)=\mathrm{I}+\sum_{m=1}^{N}\sum_{n=1}^{N}\frac{v_m(M^{-1})_{mn}\hat v_n}{k-\hat k_n}.\label{4.1b}
\end{align}
\end{subequations}
where $M=(M_{mn})_{N\times N}$ is a matrix whose entries are
\eqa M_{mn}=\frac{\hat v_mv_n}{k_m^*-k_n},\,1\leq m,n\leq N.\label{4.2}\eeqa
Therefore, from \eqref{4.1a} and \eqref{4.1b}, we obtain
\beq P_{+}^{(1)}=\sum_{m=1}^{N}\sum_{n=1}^{N}v_m(M^{-1})_{mn}\hat v_n.\label{4.3}\eeq
we have chosen $v_{n0}=[1,c_n,d_n]^T$, it follows from \eqref{4.3} that the general N-soliton solution for the mCNLS system \eqref{1.3} reads
\begin{subequations}
\begin{align}
&q_1=2\sum_{m=1}^{N}\sum_{n=1}^{N}c_me^{\theta_m-\theta_n^*}(M^{-1})_{mn},\label{4.4a}\\
&q_2=2\sum_{m=1}^{N}\sum_{n=1}^{N}d_me^{\theta_m-\theta_n^*}(M^{-1})_{mn}.\label{4.4b}
\end{align}
\end{subequations}
and $M=(M_{mn})_{N\times N}$ is given by
\eqa M_{mn}=\frac{e^{-(\theta_m^*+\theta_n)}-(c_m^*c_n-d_m^*d_n)e^{\theta_m^*+\theta_n}}{k_m^*-k_n},\,1\leq m,n\leq N.\label{4.5}\eeqa
with
$$\theta_n=-ik_nx-ik_n^2t.$$

In what follows, one can examine the nonlinear dynamical behaviors of the one-soliton solutions and two-soliton solutions to the mCNLS system \eqref{1.3}

On the one hand, as a special example, one can choose $N=1$ in formula \eqref{4.4a} and \eqref{4.4b} and with \eqref{4.2}, we obtain the one-soliton solution as follows:
\begin{subequations}
\begin{align}
&q_1(x,t)=\frac{2c_1e^{\theta_1-\theta_1^*}(k_1^*-k_1)}{e^{-(\theta_1+\theta_1^*)}-(|c_1|^2-|d_1|^2)e^{\theta_1+\theta_1^*}},\label{4.6a}\\
&q_2(x,t)=\frac{2d_1e^{\theta_1-\theta_1^*}(k_1^*-k_1)}{e^{-(\theta_1+\theta_1^*)}-(|c_1|^2-|d_1|^2)e^{\theta_1+\theta_1^*}}.\label{4.6b}
\end{align}
\end{subequations}
Letting $k_1=k_{11}+ik_{12}$, $|c_1|^2-|d_1|^2=e^{2\xi_1}$, then the one-soliton solution \eqref{4.6a} and \eqref{4.6b} can be written as
\begin{subequations}
\begin{align}
&q_1(x,t)=2ic_1k_{12}e^{\theta_1-\theta_1^*-\xi_1}\mathrm{csch}(\theta_1+\theta_1^*+\xi_1),\label{4.7a}\\
&q_2(x,t)=2id_1k_{12}e^{\theta_1-\theta_1^*-\xi_1}\mathrm{csch}(\theta_1+\theta_1^*+\xi_1).\label{4.7b}
\end{align}
\end{subequations}
We show that the one-soliton solution $q_1(x,t)$ and $q_2(x,t)$ in Fig. 1 with the parameters chose as $c=-2,d=-4,k=0.2+0.03i$.
\begin{figure}
\centering
\includegraphics[width=2.5in,height=3.5in]{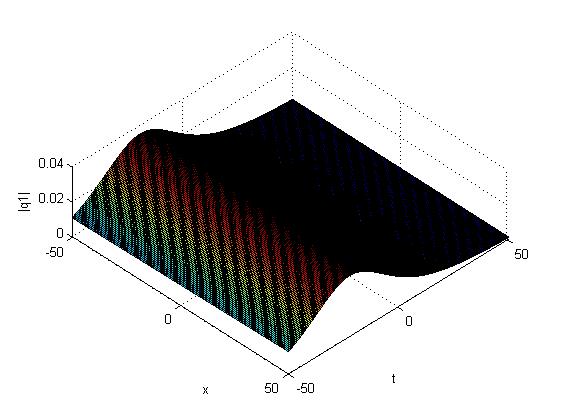}
\includegraphics[width=2.5in,height=3.5in]{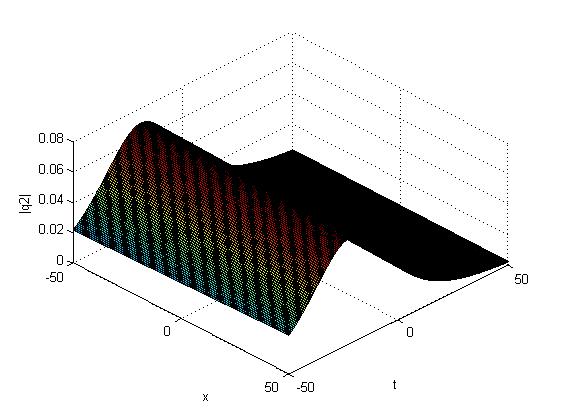}
\caption{The single bright-bright soliton solutions $q_1(x,t)$ and $q_2(x,t)$ with the parameters chose as $c=-2,d=-4,k=0.2+0.03i$.}
\label{fig:graph}
\end{figure}

On the other hand, as another special example, one can choose $N=2$ in formula \eqref{4.4a} and \eqref{4.4b} and with \eqref{4.2}, we arrive at the two-soliton solution as follows:

\begin{subequations}
\begin{align}
&q_1(x,t)=2[c_1e^{\theta_1-\theta_1^*}(M^{-1})_{11}+c_1e^{\theta_1-\theta_2^*}(M^{-1})_{12}
c_2e^{\theta_2-\theta_1^*}(M^{-1})_{21}+c_2e^{\theta_2-\theta_2^*}(M^{-1})_{22}],\label{4.8a}\\
&q_2(x,t)=2[d_1e^{\theta_1-\theta_1^*}(M^{-1})_{11}+d_1e^{\theta_1-\theta_2^*}(M^{-1})_{12}
d_2e^{\theta_2-\theta_1^*}(M^{-1})_{21}+d_2e^{\theta_2-\theta_2^*}(M^{-1})_{22}].\label{4.8b}
\end{align}
\end{subequations}
where $M=(M_{mn})_{2\times 2}$ with
$$M_{11}=\frac{e^{-(\theta_1+\theta_1^*)}-(|c_1|^2-|d_1|^2)e^{\theta_1+\theta_1^*}}{k_1^*-k_2},\,\,
M_{12}=\frac{e^{-(\theta_2+\theta_1^*)}-(c_1^*c_2-d_1^*d_2)e^{\theta_2+\theta_1^*}}{k_1^*-k_2},$$
$$M_{21}=\frac{e^{-(\theta_1+\theta_2^*)}-(c_1c_2^*-d_1d_2^*)e^{\theta_2^*+\theta_1}}{k_2^*-k_1},\,\,
M_{22}=\frac{e^{-(\theta_2+\theta_2^*)}-(|c_2|^2-|d_2|^2)e^{\theta_2+\theta_2^*}}{k_2^*-k_2},$$

\section{Discussions}

Indeed, as a promotion, the integrable two-component NLS equation or Manokov system \eqref{1.2} can be extended to the integrable generalized multi-component NLS system as follows:
\eqa i\mathbf{q}_t+\frac{1}{2}\mathbf{q}_{xx}+\epsilon \mathbf{q}\mathbf{q}^\dag \Omega\mathbf{q}=0,\epsilon=\pm1,\label{5.1}\eeqa
where $\mathbf{q}=(q_1,q_2,...,q_N)^T$,
$\Omega=diag(\omega_1,\omega_2,...,\omega_N)$,
which possess the following Lax pair for $\epsilon=-1$
\begin{subequations}
\begin{align}
&\Phi_x=(-ik\Lambda+iQ)\Phi,\label{5.2a}\\
&\Phi_t=(-ik^2\Lambda+ikQ-\frac{1}{2}(i\Lambda Q^2-\Lambda Q_x))\Phi,\label{5.2b}
\end{align}
\end{subequations}
where $k\in \mathbb{C}$ is a spectral parameter and
\eqa\begin{array}{l}
\Lambda=\left(\begin{array}{cc}
-1&\mathbf{0}_{1\times N}\\
\mathbf{0}_{N\times 1}&\mathbf{I}_{N\times N}\end{array} \right),
Q=\left(\begin{array}{ccc}
0&-\mathbf{q}^\dag S\\
\mathbf{q}&\mathbf{0}_{N\times N}\end{array} \right).
\end{array}\label{5.3}\eeqa

Indeed, if all $\omega_i=1$, which means to the focusing case, if all $\omega_i=-1$, which means to the defocusing case, or otherwise the mixed case. Accordingly, one can also examine the N-soliton solutions to the integrable generalized multi-component NLS system by the same way in Section 3. However, we don't examine them here since the procedure is mechanical.

\section{Conclusions}

In this work, we have established the multi-soliton solution of the mCNLS system on the line via the RH approach.
Basing on the Jost solutions to the Lax pair of the mCNLS system and the scattering matrix $S(k)$,
we constructed the corresponding RH problem, which admitted simple zero points
generated by the roots of det $s_{11}(k)$. By taking spectral analysis, we found that the zero points were
paired, since det$s_{11}(k)$ was an odd function. In view of the symmetry relations of zero points,
we constructed a transformation, which eliminated the zero points and made the RH
problem be regular. Applying the Plemelj formulae, N-soliton solutions to the mCNLS system were
obtained from the solutions of RH problem with vanishing scattering coefficients,
which was just the reflection-less case.
For other integrable equations with high-order matrix Lax pairs, can we construct their multi-soliton solution of associated matrix RH problem formulated in the complex $k$-plane according to RH approach?
Moreover, based on the $3\times 3$ matrix RH problem of the mCNLS system examined by Tian in \cite{TSF2016}, one can examine the long-time asymptotic behavior for the solutions of the mCNLS system through the nonlinear steepest descent method introduce by Deift and Zhou \cite{DP1993}.
This two questions will be discussed in our future paper.

\section*{Acknowledgements}

This work was supported by the National Natural Science Foundation of China under the grant No. 11601055 and 11805114, Natural Science Foundation of Anhui Province under the grant No. 1408085QA06, Education Department scientific research project of Anhui Province under the grant No. KJ2017B10.

\end{document}